\documentclass{file}
\usepackage{pat}
\usepackage[T1]{fontenc}
\usepackage[utf8]{inputenc}
\usepackage{mathtools,amssymb,tikz,float}
\usetikzlibrary{arrows.meta,shapes.geometric,decorations.pathreplacing}
\usepackage{paralist}
\usepackage[longnamesfirst,numbers,sort&compress]{natbib}
\usepackage{hyperref}
\newcommand{\edge}{\mathord{\,\,\vrule width 12pt height 3pt depth -1.5pt \,\,}}
\usepackage[mathlines]{lineno}
\newcommand*\patchAmsMathEnvironmentForLineno[1]{%
  \expandafter\let\csname old#1\expandafter\endcsname
    \csname #1\endcsname
  \expandafter\let\csname oldend#1\expandafter\endcsname
    \csname end#1\endcsname
  \renewenvironment{#1}
    {\linenomath\csname old#1\endcsname}
    {\csname oldend#1\endcsname\endlinenomath}%
}

\newcommand*\patchBothAmsMathEnvironmentsForLineno[1]{%
  \patchAmsMathEnvironmentForLineno{#1}%
  \patchAmsMathEnvironmentForLineno{#1*}%
}

\AtBeginDocument{%
  \patchBothAmsMathEnvironmentsForLineno{equation}%
  \patchBothAmsMathEnvironmentsForLineno{align}%
  \patchBothAmsMathEnvironmentsForLineno{flalign}%
  \patchBothAmsMathEnvironmentsForLineno{alignat}%
  \patchBothAmsMathEnvironmentsForLineno{gather}%
  \patchBothAmsMathEnvironmentsForLineno{multline}%
}
\newenvironment{clmproof}{\begin{proof}[Proof of Claim:]}{\end{proof}}
\newcommand{\cay}{\operatorname{Cay}}

\newcommand{\Aut}{\operatorname{Aut}}

\title{\MakeUppercase{
Hamilton cycles in vertex-transitive graphs
of prime-power order
}}
\author{
Farzad Maghsoudi,\quad  Dragan Marušič\thanks{Department of Mathematics, University of Primorska, UP IAM, Muzejski trg 2, 6000 Koper, Slovenia},\quad Bobby Miraftab
\thanks{School of Computer Science,
Carleton University, Ottawa, ON, Canada.}
}
\date{}

\begin{document}

\maketitle

\begin{abstract}
A graph is vertex-transitive if, for any two vertices, there is an automorphism mapping one to the other.
We prove that every connected vertex-transitive graph of order $p^n$, where $p$ is a prime and $p^n\geq 3$, has a hamiltonian cycle.
Our main tool extends Witte's theorem [JCTB (1986)] on Cayley digraphs of $p$-groups.
\end{abstract}

\section*{AI disclosure}

The proof presented in this paper was developed through an iterative exchange between the authors and AI.
The authors suggested developing a proof along the lines of Witte's result on Cayley digraphs of finite $p$-groups \cite{Witte1986}. 
More precisely, the assistant helped work out the details of our original idea of the cyclic-layer argument in \Cref{lem:cyclic-layers} and the connectedness argument underlying \Cref{lem:connected-lift}.
AI-generated suggestions were not accepted without verification. 
All results were checked by the authors, who take full responsibility for the content.

\section{Introduction}

A graph is \defin{vertex-transitive} if its automorphism group acts transitively on its vertices.
In 1969, Lov\'asz~\cite{Lovasz1970} asked whether every finite connected vertex-transitive graph has a hamiltonian
path. This question is still open, although it has been answered affirmatively in many special cases; see, for example,~\cite{Maghsoudi6pq,Maghsoudipqrs,LehnerMaghsoudiMiraftab2026,LehnerMaghsoudiMiraftabpqrs,MaghsoudiMiraftabp2q}.
Only four connected vertex-transitive graphs on at least three vertices are known to have no hamiltonian cycle:
the Petersen graph, the Coxeter graph, and the two graphs obtained from them by replacing each vertex with a triangle. None of these four graphs is a Cayley graph, and it is a folklore conjecture that every connected Cayley graph on at least three vertices has a hamiltonian cycle.
We refer to~\cite{KutnarMarusic2008} for a survey.

Most partial results concern graphs of restricted order. For example,
Du, Kutnar and Maru\v{s}i\v{c}~\cite{DuKutnarMarusic2021} proved that,
apart from the Petersen graph, every connected vertex-transitive graph
of order $pq$, where $p$ and $q$ are primes, has a hamiltonian cycle.
For prime-power orders, the starting point is Witte's
theorem~\cite{Witte1986} that every strongly connected Cayley digraph of
a nontrivial finite $p$-group has a directed hamiltonian cycle.
Maru\v{s}i\v{c}~\cite{Marusic1985} proved that every vertex-transitive
digraph of order $p^k$ with $k\leq 3$ is a Cayley digraph, and that
every connected vertex-transitive graph of these orders is hamiltonian.
He also constructed non-Cayley vertex-transitive graphs of order $p^k$
for every $k\geq 4$, so Witte's theorem does not settle the prime-power case. 
Chen~\cite{Chen1998} and Zhang~\cite{Zhang2015} proved that connected vertex-transitive digraphs of order $p^4$ and $p^5$, respectively, have directed hamiltonian cycles. Both proofs lift hamiltonian cycles from quotient digraphs. To the best of our knowledge, the prime-power case has remained open for every exponent $n\geq 6$.
We settle it in this paper.

\begin{thm}\label{thm:main}
Let $p$ be a prime and let $n\geq1$, with $p^n\geq3$.
Every connected vertex-transitive graph of order $p^n$ has
a hamiltonian cycle.
\end{thm}

The condition $p^n\geq3$ excludes only $K_2$. In particular, every
connected vertex-transitive graph of prime-power order has a hamiltonian
path, which answers Lov\'asz's question for these orders. Consistently
with \Cref{thm:main}, none of the four exceptional graphs above has
prime-power order; their orders are $10$, $28$, $30$ and $84$.

Our proof has two ingredients. The first concerns digraphs with cyclic
layers. Let $A$ be a finite group, let $m\geq1$, and let
$S_0,\ldots,S_{m-1}$ be nonempty subsets of $A$. The digraph
$D(A;m;S_0,\ldots,S_{m-1})$ has vertex set $A\times\mathbb Z_m$ and arcs
$(x,i)\rightarrow(xs,i+1)$ for $s\in S_i$. \Cref{lem:cyclic-layers}
shows that if $A$ is a $p$-group and this digraph is strongly connected,
then it has a directed hamiltonian cycle. The case $m=1$ is Witte's
theorem (\Cref{cor:witte-p-groups}). The proof follows Witte's strategy:
the Burnside basis theorem yields a homomorphism from $A$ onto a cyclic
group of order $p$, and this homomorphism produces a strongly connected
spanning subdigraph isomorphic to $D(N;pm;T_0,\ldots,T_{pm-1})$ for a
subgroup $N$ of index $p$ in $A$, so induction on $|A|$ applies.

The second ingredient reduces a vertex-transitive graph $X$ of order
$p^n$ to such a digraph. We choose an inclusion-minimal transitive
$p$-subgroup $P$ of $\Aut(X)$. If $P$ acts regularly, then $X$ is a
Cayley graph of $P$ and Witte's theorem applies. Otherwise, a lemma of
Zhang~\cite[Lemma~2.3]{Zhang2015} provides a central element $w$ of a
vertex stabilizer and a nontrivial normal subgroup $N$ of $P$ that acts
semiregularly on $V(X)$. The quotient of $X$ by the $N$-orbits is a
connected vertex-transitive graph of smaller prime-power order. By
induction it has a hamiltonian cycle $B_0B_1\cdots B_{m-1}$, except when
$m=2$, which we handle directly. The arcs of $X$ from $B_i$ to $B_{i+1}$
form a spanning digraph isomorphic to $D(N;m;S_0,\ldots,S_{m-1})$. The
main step, \Cref{lem:connected-lift}, uses $w$ to show that this digraph
is strongly connected, and \Cref{lem:cyclic-layers} completes the proof.

Section~2 treats digraphs with cyclic layers and recovers Witte's theorem. Section~3 proves \Cref{thm:main}.

All graphs are finite and simple, and undirected unless we call them digraphs. A \defin{hamiltonian cycle} visits every vertex exactly once
before returning to its starting vertex; directed hamiltonian cycles of digraphs are defined analogously. A digraph is \defin{strongly connected}
if every vertex can reach every other vertex by a directed walk.

\section{Results}
\begin{defn}
Let $A$ be a finite group, let $m\geq 1$, and let
$\varnothing\neq S_i\subseteq A$ for each $i\in\mathbb Z_m$.
The digraph $D=D(A;m;S_0,\ldots,S_{m-1})$
has vertex set $A\times\mathbb Z_m$ and arcs
$(x,i)\rightarrow(xs,i+1)$, where 
$x\in A,\ i\in\mathbb Z_m,\ s\in S_i$,
and $i+1$ is modulo $m$.
We call $A\times\{i\}$ the \emph{$i$th layer} and label the
displayed arc by $s$.
Thus an arc multiplies the group coordinate on the right by $s$
and advances to the next layer; layer $m-1$ leads back to layer $0$.
\end{defn}

\begin{exa}
Take $A=\Z_2=\{e,a\}, m=3,S_0=\{e,a\}, S_1=S_2=\{e\}$.
There are three layers, each containing two vertices.
From $(x,0)$ there are arcs to both $(x,1)$ and $(xa,1)$.
The remaining arcs are
\[
 (x,1)\rightarrow(x,2)\rightarrow(x,0)
 \qquad(x\in\{e,a\}).
\]
An arc labelled $e$ preserves the group coordinate.
An arc labelled $a$ switches $e$ and $a$, since $ea=a$ and $aa=e$, see \Cref{fig:exa}.
\begin{figure}[H]
\centering
\includegraphics[scale=0.9]{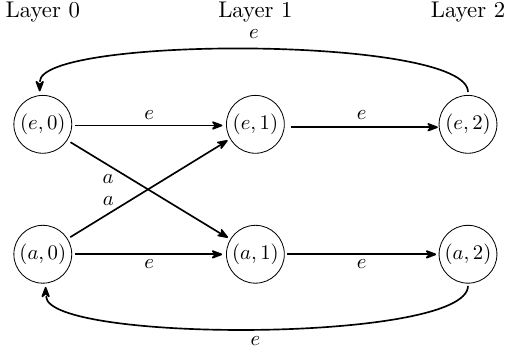}
\caption{All eight arcs are shown. The crossing of the diagonal arcs is not a vertex.}\label{fig:exa}
\end{figure}
\end{exa}

Next we use Burnside basis theorem.
\begin{lem}\label{lem:burnside}
Let $P$ be a finite $p$-group and let $\Phi(P)$ be its
Frattini subgroup. Then $P/\Phi(P)$ is an $\mathbb F_p$-vector
space, and $T\subseteq P$ generates $P$ if and only if its image spans $P/\Phi(P)$.
\end{lem}

\begin{lem}\label{lem:cyclic-layers}
Let $A$ be a finite $p$-group, $m\geq 1$, and
$\varnothing\neq S_i\subseteq A$ for $i\in\mathbb Z_m$.
The digraph $D=D(A;m;S_0,\ldots,S_{m-1})$
has a directed hamiltonian cycle whenever it is strongly connected.
For the one-vertex case, the available loop is considered as the cycle.
\end{lem}

\begin{proof}
We induct on $|A|$, for all $m$ simultaneously.
If $|A|=1$, the layers themselves form the required directed
hamiltonian cycle, including the loop when $m=1$.

Put $R=S_0\cdots S_{m-1}$. 
These are precisely the words making one complete trip through the layers.
We first claim the following:
\begin{clm}\label{clm:scc}
$D$    is strongly connected if and only if $\langle R\rangle=A$.
\end{clm}

\begin{clmproof}
Each product $r=s_0s_1\cdots s_{m-1}\in R$ determines
a directed walk of length $m$:
\[
(x,0)\rightarrow(xs_0,1)
\rightarrow\cdots\rightarrow(xr,0).
\]

We first suppose that $D$ is strongly connected.
For every $a\in A$, there is a directed walk from $(e,0)$
to $(a,0)$, where $e$ is the identity of $A$.
Every arc increases the layer index by $1$ modulo $m$,
so this walk has length $km$ for some $k\geq 0$.

Split the walk into $k$ consecutive blocks of length $m$.
Each block starts in layer $0$ and uses, in order,
one label from each of $S_0,\ldots,S_{m-1}$.
Its product therefore belongs to $R$.
Writing these products as $r_1,\ldots,r_k$, we obtain
$a=r_1r_2\cdots r_k\in\langle R\rangle$.
Since $a$ was arbitrary, $\langle R\rangle=A$.

Conversely assume that $\langle R\rangle=A$.
For any $x,y\in A$, we can write $x^{-1}y=r_1r_2\cdots r_k$, where $r_1,\ldots,r_k\in R$.
Concatenating the directed walks of length $m$ corresponding
to $r_1,\ldots,r_k$ gives
\[
(x,0)\leadsto(xr_1,0)
\leadsto\cdots\leadsto(xr_1\cdots r_k,0)=(y,0).
\]
Therefore every vertex in layer $0$ can reach every other
vertex in layer $0$.

Now consider arbitrary vertices $(x,i)$ and $(y,j)$,
with $0\leq i,j<m$.
Choose $a_\ell\in S_\ell$ for $0\leq\ell<m$, and put
\[
u=
\begin{cases}
e, & i=0,\\
a_i a_{i+1}\cdots a_{m-1}, & i>0,
\end{cases}
\qquad
v=
\begin{cases}
e, & j=0,\\
a_0a_1\cdots a_{j-1}, & j>0.
\end{cases}
\]
Following the labels $a_i,\ldots,a_{m-1}$ when $i>0$
gives a directed walk $(x,i)\leadsto(xu,0)$.
By the preceding paragraph, there is a directed walk
$(xu,0)\leadsto(yv^{-1},0)$.
Finally, following $a_0,\ldots,a_{j-1}$ when $j>0$ gives
$(yv^{-1},0)\leadsto(yv^{-1}v,j)=(y,j)$.
When $i=0$ or $j=0$, the corresponding walk has length zero.
Concatenating these walks shows that $(x,i)$ can reach
$(y,j)$. Hence $D$ is strongly connected.
\end{clmproof}

Write $A/\Phi(A)$ multiplicatively, and denote the image
of $x\in A$ by $\bar x$. By Lemma~\ref{lem:burnside},
$A/\Phi(A)$ is elementary abelian and $\bar R$ generates $A/\Phi(A)$.
Let $b\in R\sm \Phi(A)
$.
Choose $a_i\in S_i$ such that $b=a_0\cdots a_{m-1}$.
Thus $\bar b\neq1$. Since $A/\Phi(A)$ is abelian, for every
choice $s_i\in S_i$,
\[
\overline{s_0\cdots s_{m-1}}
 =\prod_{i=0}^{m-1}\bar s_i
 =\bar b\prod_{i=0}^{m-1}
       \bigl(\bar s_i\bar a_i^{-1}\bigr).
\]
Consequently, we have $\{\bar b\}\cup \{\bar s\bar a_i^{-1}
\colon 0\leq i<m,\ s\in S_i\}$ generate $A/\Phi(A)$.

Choose among these elements $u_1,\ldots,u_k$ so that
$\bar b,u_1,\ldots,u_k$ form a minimal generating set
of $A/\Phi(A)$, allowing $k=0$.
Since $A/\Phi(A)$ is elementary abelian, every element of $A/\Phi(A)$
has a unique expression $\bar b^{\,t_0}u_1^{t_1}\cdots u_k^{t_k}$, where $0\leq t_0,\ldots,t_k<p$.
Let $C_p=\langle\zeta\rangle$ be a cyclic group of order $p$. We define the following homomorphism
\begin{align*}
\lambda:\quad A/\Phi(A) &\rightarrow C_p,\\
\bar b^{\,t_0}u_1^{t_1}\cdots u_k^{t_k}
&\longmapsto \zeta^{t_0}.
\end{align*}
We notice that $\lambda(\bar b)=\zeta, \lambda(u_j)=1$.
Next we define the homomorphism $\phi(x)=\lambda(\bar x)$,
Assume that 
$ c_i=\phi(a_i),
S_i'=\{s\in S_i:\phi(s)=c_i\}$
and let $D'=D(A;m;S_0',\ldots,S_{m-1}')$.

Each $S_i'$ contains $a_i$.
Moreover, if a selected element has the form
$u_j=\bar s\bar a_i^{-1}$, then
$1=\lambda(u_j)=\phi(s)\phi(a_i)^{-1}$,
so $s\in S_i'$.
Next we set $R'=S_0'\cdots S_{m-1}'$.
Since $a_i\in S_i'$ for every $i$, we have
$b=a_0\cdots a_{m-1}\in R'$.
For each selected element $u_j=\bar s\bar a_i^{-1}$,
we have already shown that $s\in S_i'$. Hence
$a_0\cdots a_{i-1}s\,a_{i+1}\cdots a_{m-1}\in R'$.
Since $A/\Phi(A)$ is abelian, the image of this product is
\[
\overline{a_0\cdots a_{i-1}s\,a_{i+1}\cdots a_{m-1}}
 =\bar b\,\bar s\bar a_i^{-1}
 =\bar b\,u_j.
\]
Thus $\overline{R'}$ contains $\bar b$ and every
$\bar b u_j$. These elements generate $A/\Phi(A)$, because
$u_j=\bar b^{-1}(\bar b u_j)$
and $\bar b,u_1,\ldots,u_k$ generate $A/\Phi(A)$.
Therefore \Cref{lem:burnside} gives $\langle R'\rangle=A$,
and \Cref{clm:scc} shows that $D'$ is strongly connected.
Since $\phi(b)=\zeta$, the homomorphism $\phi$ is
surjective. Consequently, we obtain 
$|N|=|A|/p$, where $N=\ker\phi$.
For $z\in C_p$ and $i\in\mathbb Z_m$, define
\[
B_{z,i}=\{(x,i):\phi(x)=z\}.
\]
These sets partition the vertices of $D'$.
If $(x,i)\in B_{z,i}$ and $s\in S_i'$, then
$\phi(xs)=\phi(x)\phi(s)=zc_i$.
Thus every arc of $D'$ starting in $B_{z,i}$ ends in
$B_{zc_i,i+1}$. Since each $S_i'$ is nonempty, the quotient digraph has precisely the arcs $B_{z,i}\rightarrow B_{zc_i,\,i+1}$.
Furthermore, we have 
\[
c_0c_1\cdots c_{m-1}
 =\phi(a_0\cdots a_{m-1})
 =\phi(b)=\zeta.
\]
Starting at $B_{1,0}$, after $qm$ steps the quotient walk
therefore reaches $B_{\zeta^q,0}$.
A return to $B_{1,0}$ requires the number of steps to be
a multiple of $m$, and then requires $\zeta^q=1$.

Since $\zeta$ has order $p$, the first positive return
occurs after $pm$ steps.
There are exactly $pm$ sets $B_{z,i}$, so the quotient
is one directed cycle through all of them, see \Cref{fig:placeholder}.
\begin{figure}[H]
\centering\includegraphics[width=0.5\linewidth]{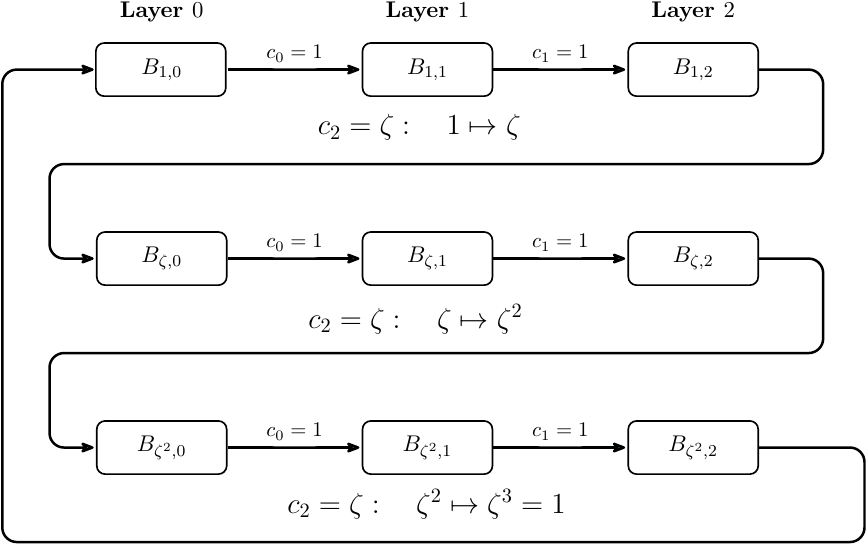}
    \caption{Each trip takes $m$ steps and multiplies $z$ by
   $c_0\cdots c_{m-1}=\zeta$.
   Since $\zeta$ has order $p$, the first return takes $p$ trips and so totally$pm$ steps. In this example, we have $p=m=3$, $c_0=c_1=1$, $c_2=\zeta$, $|\zeta|=3$.}
    \label{fig:placeholder}
\end{figure}
Put $\ell=pm$, and list these sets in their directed order
as $B_0,\ldots,B_{\ell-1}$.
If $B_j=B_{z_j,i_j}$, choose $r_j\in A$ with
$\phi(r_j)=z_j$. Since $N=\ker\phi$,
we have $B_j=\{(nr_j,i_j):n\in N\}$.
We consider $j$ modulo $\ell$, so $r_\ell=r_0$.
Next we define $T_j=r_jS_{i_j}'r_{j+1}^{-1}$.
Each $T_j$ is nonempty.
For every $s\in S_{i_j}'$, the arc from $(r_j,i_j)$
labelled $s$ ends in $B_{j+1}$. Hence we have $r_js\in Nr_{j+1}$, and so
$r_jsr_{j+1}^{-1}\in N$.
Therefore we obtain $T_j\subseteq N$.
Now we consider $E=D(N;\ell;T_0,\ldots,T_{\ell-1})$.
We now define
\begin{align*}
\Psi\colon V(E)&\rightarrow V(D'),\\
(n,j)&\longmapsto(nr_j,i_j).
\end{align*}
The sets $B_j$ partition $V(D')$, and each vertex of
$B_j$ has a unique expression $(nr_j,i_j)$.
Thus $\Psi$ is a bijection.

For $s\in S_{i_j}'$, put
$t=r_jsr_{j+1}^{-1}\in T_j$. Then
we have $(nr_j)s=(nt)r_{j+1}$.
Consequently, the arc
$(n,j)\rightarrow(nt,j+1)$
of $E$ corresponds under $\Psi$ to the arc $(nr_j,i_j)\rightarrow((nr_j)s,i_{j+1})$
of $D'$. Every element of $T_j$ arises from some
$s\in S_{i_j}'$, so this correspondence preserves arcs
in both directions. Hence $E$ and $D'$ are isomorphic.

In particular, $E$ is strongly connected.
Since $|N|=|A|/p<|A|$, induction gives a directed
hamiltonian cycle in $E$.
Its image under $\Psi$ is a directed hamiltonian cycle
in $D'$, and therefore also in $D$.
\end{proof}

\begin{cor}{\rm\cite{Witte1986}
\label{cor:witte-p-groups}}
Every strongly connected Cayley digraph of a non-trivial finite $p$-group has a directed hamiltonian cycle.
\end{cor}

\begin{proof}
We invoke \Cref{lem:cyclic-layers} with $m=1$ and $S_0=S$.
Under the identification $(x,0)\mapsto x$, the digraph
$D(A;1;S)$ is precisely
$\cay(A;S)$.
\end{proof}

%%%%%%%%%%%%%%%%%%%%%%%%%%%%%%%%%%%%%%%%%%%%%%%%%%%%%%%%%%%%%%%%%%%%%%%%%%%%%%%%%%%%%%%%%%%%%

\section{Vertex-transitive graphs}

All groups in this section act on the left. If a group $P$ acts
transitively on a set and $H$ is the stabilizer of a point, we identify
the set with the left cosets $P/H$. Thus $a\in P$ sends $gH$ to $agH$.
We first recall three group-theoretic results.

\begin{lem}[{\cite[Theorems~5.2 and~5.6(c), pp.~78--79]{Milne}}]
\label{lem:sylow}
Let $G$ be a finite group and let $p$ be a prime. Every $p$-subgroup
of $G$ is contained in a subgroup of order $|G|_p$, where $|G|_p$
is the largest power of $p$ dividing $|G|$.
\end{lem}

\begin{lem}[{\cite[\S6, p.~14]{Conrad}}]\label{lem:maximal}
Every maximal proper subgroup of a finite $p$-group is normal and
has index $p$.
\end{lem}

\begin{lem}[{\cite[Lemma~2.3, pp.~3--4]{Zhang2015}}]\label{lem:zhang}
Let $P$ be a finite $p$-group and let $\{1\}<H<P$, with
$H\cap Z(P)=\{1\}$. There is an element $w\in Z(H)$ of order $p$ such that $N=\langle wgw^{-1}g^{-1}:g\in P\rangle$
satisfies  $N\cap H=\{1\}$.
\end{lem}

%Recall that $[g,u]=g^{-1}u^{-1}gu$.
%If $u$ is the element supplied by his lemma, take $w=u^{-1}$.
%As $g$ ranges over $P$, the elements $wgw^{-1}g^{-1}$ are the
%inverses of the elements $[g^{-1},u]$, so they generate the same
%subgroup. We have stated explicitly the necessary hypothesis
%$H\ne\{1\}$.

\begin{lem}\label{lem:transitive-p}
A finite transitive permutation group of degree $p^n$ contains a transitive $p$-subgroup. Let $P$ be an inclusion-minimal transitive
$p$-subgroup, and let $H$ be a point stabilizer in $P$. Then
\[
H\le\Phi(P)\qquad\text{and}\qquad H\cap Z(P)=\{1\}.
\]
\end{lem}

\begin{proof}
Let $G$ be the given permutation group and let $v$ be a point.
By \Cref{lem:sylow}, choose $S\le G_v$ of order $|G_v|_p$,
and then choose $P\le G$ of order $|G|_p$ with $S\le P$.
Since $P\cap G_v$ is a $p$-subgroup of $G_v$ containing $S$,
we have $P\cap G_v=S$. Consequently,
\[
|P:P\cap G_v|
=\frac{|G|_p}{|G_v|_p}
=|G:G_v|=p^n.
\]
Thus the orbit of $v$ under $P$ contains every point.

Now choose $P$ inclusion-minimal among the transitive $p$-subgroups,
and put $H=P_v$. Suppose a maximal proper subgroup $M<P$ does
not contain $H$. By \Cref{lem:maximal}, $M$ is normal and has
index $p$, so $MH=P$. Every coset $gH$ therefore has a
representative in $M$, making $M$ transitive on $P/H$.
This contradicts the choice of $P$. Hence $H$ lies in every
maximal proper subgroup of $P$, and $H\le\Phi(P)$.

Finally, if $h\in H\cap Z(P)$, then
$h(gv)=g(hv)=gv$ for every $g\in P$. Thus $h$ fixes every point
and, as $P$ is a permutation group, $h=1$.
\end{proof}

\begin{lem}\label{lem:cycle-components}
Let $Y$ be a finite digraph and let $N\leq\Aut(Y)$.
Suppose that only the identity of $N$ fixes a vertex.
Let $B_0,\ldots,B_{m-1}$ be the $N$-orbits, where $m\geq2$.
Assume that every arc goes from some $B_i$ to $B_{i+1}$,
and that, for every $i\in\mathbb Z_m$, there exist
$x\in B_i$ and $y\in B_{i+1}$ such that
$x\rightarrow y$ is an arc of $Y$.
Then every component of the underlying graph
induces a strongly connected subdigraph of $Y$ and meets
every $B_i$.
\end{lem}

\begin{proof}
For an arc $x\rightarrow y$ from $B_i$ to $B_{i+1}$,
the arcs $nx\rightarrow ny$ for $n\in N$
form a perfect matching between these two sets.
%Indeed, the maps $n\mapsto nx$ and $n\mapsto ny$ are
%bijections onto $B_i$ and $B_{i+1}$, respectively.
Fix any arc of $Y$. 
Suppose the fixed arc goes from $B_j$ to $B_{j+1}$.
Choose the perfect matching containing this arc.
For each $i\in\Z_m$ with $i\neq j$, choose one perfect
matching from $B_i$ to $B_{i+1}$.
Their union has exactly one incoming and one outgoing
arc at each vertex, so it is a disjoint union of directed
cycles. Thus every arc of $Y$ lies on a directed cycle.

Consequently, the endpoints of every arc can reach each
other. Applying this along an undirected path shows that
each component is strongly connected.
Finally, every vertex has an outgoing arc to the next $B_i$. 
Following these arcs shows that each component
meets every $B_i$.
\end{proof}

\begin{lem}\label{lem:component-criterion}
Let $P$ be a finite $p$-group acting transitively on a
finite set $\Omega$. Fix $v\in\Omega$ and choose
$w\in Z(P_v)$. 
Let $N=\langle wgw^{-1}g^{-1} :g\in P\rangle$.
Let $F$ be a graph on $\Omega$ such that $N\cup \{w\}\subseteq \Aut (F)$.
If every component of $F$ meets every $N$-orbit,
then $F$ is connected.
\end{lem}

\begin{proof}
Let $C$ be the vertex set of the component containing $v$.
Since $w$ fixes $v$, it preserves $C$.
We put $J=\langle N,w\rangle$, and $K=\{j\in J:jC=C\}$.
Since $N\cup \{w\}\subseteq \Aut (F)$, $K$ is a subgroup of $J$, and $w\in K$.
Suppose that $C\neq\Omega$.
If $K=J$, every element of $N$
would preserve $C$. Since $C$ meets every $N$-orbit,
it would contain every such orbit, giving $C=\Omega$.
So we have shown that $K<J$.
Choose a maximal proper subgroup $M<J$ containing $K$.
By \Cref{lem:maximal}, we infer that $M\trianglelefteq J$.

Choose one vertex $v_i\in C$ from each $N$-orbit,
and choose $g_i\in P$ such that $g_iv=v_i$.
Since $w$ fixes $v$, we have $wg_iw^{-1}g_i^{-1}v_i=wv_i\in C$.
The automorphism $wg_iw^{-1}g_i^{-1}\in N$ therefore sends
a vertex of $C$ to another vertex of $C$.
It follows that $wg_iw^{-1}g_i^{-1}$ preserves this component.
Hence $wg_iw^{-1}g_i^{-1}\in K\leq M$.

Now let $g\in P$. The vertex $gv$ belongs to some orbit
$Nv_i$, so $gv=ng_iv$ for some $n\in N$.
Therefore $g=ng_i h$ for some $h\in P_v$.
Since $w\in Z(P_v)$, it commutes with $h$.
Direct multiplication gives
\begin{equation}\label{eq:delta-decomposition}
wgw^{-1}g^{-1}
=(wnw^{-1}n^{-1})
 \bigl(nwg_iw^{-1}g_i^{-1}n^{-1}\bigr).
\end{equation}
Both factors belong to $M$: the first because $w\in M$
and $M\trianglelefteq J$, and the second because
$wg_iw^{-1}g_i^{-1}\in M$ and $n\in J$.
Thus $wgw^{-1}g^{-1}\in M$ for every $g\in P$, giving
$N=\langle wgw^{-1}g^{-1}:g\in P\rangle\leq M$.
Since $w\in M$ as well, we obtain
$J=\langle N,w\rangle\leq M$, contrary to $M<J$.
Therefore $C=\Omega$, and $F$ is connected.
\end{proof}

\begin{lem}\label{lem:connected-lift}
Let $X$ be a graph, and let
$P\leq\Aut(X)$ be a transitive $p$-group.
Fix a vertex $v$ and choose $w\in Z(P_v)$ of order $p$.
Let
$N=\langle wgw^{-1}g^{-1}:g\in P\rangle$.
Assume that $N\cap P_v=\{1\}$.
Suppose the $N$-orbits are
$B_0,\ldots,B_{m-1}$, where $m\geq2$, so that an edge
joins $B_i$ to $B_{i+1}$ for every $i\in\Z_m$.
Replace every edge of $X$ by its two opposite arcs,
and keep all arcs from $B_i$ to $B_{i+1}$.
Then the resulting digraph $Y$ is strongly connected.
\end{lem}

\begin{proof}
Put $\delta(g)=wgw^{-1}g^{-1}$.
For $a,g\in P$, we have
\begin{equation}\label{eq:normality}
a\delta(g)a^{-1}=\delta(a)^{-1}\delta(ag).
\end{equation}
Thus $N$ is normal in $P$.
Suppose $n\in N$ fixes a vertex $gv$.
Then $g^{-1}ng$ fixes $v$, so
$g^{-1}ng\in N\cap P_v=\{1\}$.
Hence $n=1$.

The group $N$ preserves $Y$.
By \Cref{lem:cycle-components}, every component of its
underlying undirected graph is strongly connected in $Y$
and meets every $N$-orbit.

Moreover, since $w$ fixes $v$,
$w(gv)=\delta(g)(gv)$ for $g\in P$.
Thus $w$ preserves every $N$-orbit.
As $w$ is an automorphism of $X$, it also preserves $Y$.

Apply \Cref{lem:component-criterion} to the underlying
undirected graph of $Y$. It is connected.
Since each of its components is strongly connected in $Y$,
the digraph $Y$ is strongly connected.
\end{proof}

\begin{lem}[Sabidussi {\cite[Lemma~4]{Sabidussi1958}}]
\label{lem:sabidussi}
Let $X$ be a finite simple graph, and let
$G\leq\Aut(X)$ act transitively on $V(X)$.
If $G_v=\{1\}$ for some vertex $v$, then
$X\cong\operatorname{Cay}(G;S)$, where $S=\{s\in G:v\sim sv\}$
\end{lem}

\begin{proof}[Proof of \Cref{thm:main}]
Let $X$ be a connected vertex-transitive graph of order $p^n\geq3$.
Fix $p$ and induct on $n$, assuming the result for all smaller
powers of $p$ that are at least three.
By \Cref{lem:transitive-p}, choose an inclusion-minimal
transitive $p$-subgroup $P\le\Aut(X)$.
Fix a vertex $v$.
Then $P_v\le\Phi(P)$ and $P_v\cap Z(P)=\{1\}$.

If $P_v=\{1\}$, then \Cref{lem:sabidussi} gives
$X\cong\operatorname{Cay}(P,S)$ for some
$S=S^{-1}\subseteq P\setminus\{1\}$.
Replacing every edge by its two opposite arcs therefore
gives a digraph isomorphic to $D(P;1;S)$.
This digraph is strongly connected because $X$ is connected.
By \Cref{lem:cyclic-layers}, it has a directed hamiltonian
cycle. Since $|V(X)|\ge3$, this gives a hamiltonian cycle of $X$.

We may therefore assume $P_v\ne\{1\}$.
We invoke \Cref{lem:zhang} with the subgroup $P_v$ to choose
$w\in Z(P_v)$ of order $p$ and $N=\langle wgw^{-1}g^{-1}:g\in P\rangle$ with $N\cap P_v=\{1\}$.
Equation~\eqref{eq:normality} shows that $N$ is normal in $P$.
Also, $N\ne\{1\}$, since otherwise $w$ would commute with
every element of $P$, contrary to $P_v\cap Z(P)=\{1\}$.
By \Cref{lem:burnside}, the group $P/\Phi(P)$ is abelian.
Every generator $wgw^{-1}g^{-1}$ of $N$ therefore belongs to $\Phi(P)$, so we have $K=NP_v\le\Phi(P)<P$.
As shown in the proof of \Cref{lem:connected-lift}, no non-identity element of $N$ fixes a vertex.
Write $|N|=p^a$, where $a\ge1$.
Every $N$-orbit has $p^a$ vertices, and the number of these
orbits is
\[
m=|P:K|=p^{n-a}\ge p.
\]
In particular, $m<p^n$ and $m\ne1$.

Let $Q$ be the simple graph whose vertices are the $N$-orbits,
with two distinct orbits adjacent when an edge of $X$ joins
them. 
Projecting walks of $X$ onto the orbits, and deleting
consecutive repetitions, shows that $Q$ is connected.
Since $N$ is normal in $P$, the action of $P$ on $X$ induces a transitive action on $Q$. 
Thus $Q$ is vertex-transitive.
If $m\ge3$, the induction hypothesis gives a hamiltonian cycle in $Q$; list its vertices in order as $B_0,\ldots,B_{m-1}$.
If $m=2$, then $Q$ consists of a single edge, and we use the two opposite arcs $B_0\rightarrow B_1$ and $B_1\rightarrow B_0$.

Replace every edge of $X$ by its two opposite arcs, and keep all arcs from $B_i$ to $B_{i+1}$ for $i\in\Z_m$.
By \Cref{lem:connected-lift}, the resulting spanning digraph $Y$ is strongly connected.
We choose a vertex $v_i\in B_i$ for each $i$, with indices read modulo $m$, and put
\[
S_i=\{s\in N:(v_i,sv_{i+1})\text{ is an arc of }Y\}.
\]
Each $S_i$ is nonempty, because every vertex of $B_i$ has an outgoing arc to $B_{i+1}$.
Every vertex of $B_i$ has a unique expression $nv_i$ with
$n\in N$. Hence the map
\begin{align*}
\Psi\colon N\times\Z_m&\rightarrow V(Y),\\
(n,i)&\longmapsto nv_i
\end{align*}
is a bijection.
Since $N$ preserves $Y$, its arcs from $B_i$ to $B_{i+1}$ are exactly $nv_i\rightarrow ns v_{i+1}$, where $n\in N,\ s\in S_i$.
Indeed, translating an arbitrary arc from $nv_i$ by $n^{-1}$ gives an arc from $v_i$, and its other endpoint determines the unique element $s\in S_i$.
Thus $\Psi$ is an isomorphism from $D(N;m;S_0,\ldots,S_{m-1})$ to $Y$.

By \Cref{lem:cyclic-layers}, $Y$ has a hamiltonian
dicycle. 
Its vertices, in their directed order, visit every
vertex of $X$ exactly once before returning to the first.
Every consecutive pair, including the last and first, is
joined by an edge of $X$. 
Since $p^n\ge3$, this ordering is a hamiltonian cycle of the simple undirected graph $X$.
This completes the induction, with the quotient of order
two handled directly above.
\end{proof}

%%%%%%%%%%%%%%%%%%%%%%%%%%%%%%%%%%%%%%%%%%%%%%%%%%%%%%%%%%%%%%%%%%%%%%%%%%%%%%

% Add \usetikzlibrary{arrows.meta,decorations.pathreplacing} to your preamble.

The following example makes the lifting argument explicit for a non-Cayley vertex-transitive graph of order $p^4$.

We describe the semiregular subgroup, a hamiltonian
cycle in the quotient, and the label sets of the resulting spanning lift.

\begin{exa}\label{exa:coset-p4}
Let $p\ge5$ be prime, and let $G=\langle a,b\rangle$ be a group of order $p^5$ and exponent $p$. 
We set $c=[a,b], d=[a,c], e=[b,c]$,
where $[x,y]=x^{-1}y^{-1}xy$, and assume that $d,e\in Z(G)$.
Put $H=\langle c\rangle$, and define a graph $X$ on $G/H$ by
\[
gH\edge kH
\quad\Longleftrightarrow\quad
g^{-1}k\in
H\{a,a^{-1},b,b^{-1},ab,(ab)^{-1}\}H.
\]
We apply the cyclic-layer argument to this graph, see \Cref{fig:coset-p4}.

\begin{figure}[H]
\centering
% Requires tikz, graphicx, and the TikZ libraries
% arrows.meta and decorations.pathreplacing.
\begin{minipage}[t]{0.49\linewidth}
\centering
\vspace{0pt}
\resizebox{\linewidth}{!}{%
\begin{tikzpicture}[
  x=1.05cm,y=1cm,font=\small,
  vertex/.style={circle,fill=black,inner sep=1.8pt},
  qcycle/.style={-{Stealth[length=2mm]},red!80!black,thick}
]
  % The vertex in position r of row j is (-j+r,j) in F_5^2.
  \foreach \j in {0,...,4}{
    \foreach \r in {0,...,4}{
      \pgfmathtruncatemacro{\i}{mod(\r-\j+5,5)}
      \node[vertex,label=above:{$(\i,\j)$}]
        (v\j-\r) at (\r,{-1.12*\j}) {};
    }
    \foreach \r in {0,...,3}{
      \pgfmathtruncatemacro{\s}{\r+1}
      \draw[qcycle] (v\j-\r) -- (v\j-\s);
    }
  }
  % Each row ends with a (0,1)-step to the next row.
  \foreach \j in {0,...,3}{
    \pgfmathtruncatemacro{\k}{\j+1}
    \draw[qcycle,rounded corners=4pt]
      (v\j-4) -- (4.45,{-1.12*\j})
      -- (4.45,{-1.12*\j-0.52})
      -- (-0.45,{-1.12*\j-0.52})
      -- (-0.45,{-1.12*\k}) -- (v\k-0);
  }
  % The final (0,1)-step closes the cycle.
  \draw[qcycle,rounded corners=5pt]
    (v4-4) -- (4.88,-4.48) -- (4.88,0.68)
    -- (-0.88,0.68) -- (-0.88,0) -- (v0-0);
\end{tikzpicture}%
}
\par\smallskip
{\small (a) The chosen quotient Hamiltonian cycle for $p=5$.}
\end{minipage}\hfill
\begin{minipage}[t]{0.48\linewidth}
\centering
\vspace{0pt}
\resizebox{\linewidth}{!}{%
\begin{tikzpicture}[
  font=\small,
  layer/.style={draw,rounded corners=3pt,
    minimum width=0.95cm,minimum height=0.8cm},
  layerarc/.style={-{Stealth[length=2mm]},thick}
]
  \node[layer,label=above:{$B_{jp}$}] (A) at (0,0) {$N$};
  \node[layer,label=above:{$B_{jp+1}$}] (B) at (1.55,0) {$N$};
  \node (dots) at (2.8,0) {$\cdots$};
  \node[layer,label=above:{$B_{jp+p-1}$}] (C) at (4.05,0) {$N$};
  \node[layer,label=above:{$B_{(j+1)p}$}] (D) at (6.1,0) {$N$};
  \draw[layerarc] (A) -- node[below] {$\langle d\rangle$} (B);
  \draw[layerarc] (B) -- node[below] {$\langle d\rangle$} (dots);
  \draw[layerarc] (dots) -- node[below] {$\langle d\rangle$} (C);
  \draw[layerarc] (C) -- node[below] {$\langle e\rangle$} (D);
  % B_{(j+1)p} is outside the brace: it belongs to the next row.
  \draw[decorate,decoration={brace,mirror,amplitude=4pt}]
    (-0.475,-0.7) -- (4.525,-0.7)
    node[midway,below=7pt] {row $j$: $p$ layers};
\end{tikzpicture}%
}
\par\smallskip
{\small (b) The spanning lift
$Y\cong D(N;p^2;S_0,\ldots,S_{p^2-1})$.}
\end{minipage}

\caption{The quotient cycle and its spanning lift in
\Cref{exa:coset-p4}. In (a), coordinates lie in $\mathbb F_5$, and only the edges of the chosen cycle are shown. In (b), each layer is identified with $N\cong C_p^2$; arrows specify the permitted labels between consecutive layers, and subscripts are  modulo $p^2$. The brace covers exactly the $p$ layers of row $j$.}
\label{fig:coset-p4}
\end{figure}
The commutator relations allow every element of $G$ to be written as $a^i b^j c^t d^k e^\ell$, where $i,j,t,k,\ell\in\mathbb F_p$.
Since $|G|=p^5$, these expressions are unique.
In particular, $|H|=p$, and every coset has a unique expression $a^i b^j d^k e^\ell H$. Thus $X$ has $p^4$ vertices.
Left multiplication by $G$ acts transitively on $X$, and $X$ is
connected because $a,b$ generate $G$.
Identify the coset $a^i b^j d^k e^\ell H$ with
\[
(i,j,k,L),\qquad L=\ell-i\binom j2,
\]
where $\binom j2=j(j-1)/2$ in $\mathbb F_p$.
All four coordinates are modulo $p$.
The identities
\[
b^ja=ab^jc^{-j}e^{\binom j2},
\qquad
c^ta=ac^td^{-t},
\qquad
c^tb=bc^te^{-t}
\]
follow from the commutator relations, the first by induction on $j$.
Collecting $gc^qaH$, $gc^qbH$, and $gc^qabH$, where
$g=a^ib^jd^ke^\ell$ and $q\in\mathbb F_p$, gives the following adjacency relations, together with their reverses:
\begin{equation}\label{eq:coset-p4-edges}
\begin{aligned}
(i,j,k,L)&\edge(i+1,j,k+t,L),\\
(i,j,k,L)&\edge(i,j+1,k,L+t),\\
(i,j,k,L)&\edge(i+1,j+1,k+t,L+t-ij),
\end{aligned}
\qquad t\in\mathbb F_p.
\end{equation}

We choose $w=c\in Z(H)$ in the lifting argument, and set $N=\langle cgc^{-1}g^{-1}:g\in G\rangle$.
Since $cac^{-1}a^{-1}=d^{-1}$, and $cbc^{-1}b^{-1}=e^{-1}$,
we have $\langle d,e\rangle\le N$.
Conversely, the image of $c$ is central in $G/\langle d,e\rangle$, since it commutes with the images of both generators $a,b$.
Thus $cgc^{-1}g^{-1}\in\langle d,e\rangle$ for every $g\in G$.
Consequently, we obtain  $N=\langle d,e\rangle\cong C_p^2$.
Uniqueness of the normal forms gives $N\cap H=\{1\}$.
Moreover, $d^ue^v$ acts by $(i,j,k,L)\longmapsto(i,j,k+u,L+v)$.
Hence $N$ acts semiregularly, with orbits
$B_{i,j}=\{(i,j,k,L):k,L\in\mathbb F_p\}$.
By \eqref{eq:coset-p4-edges}, the quotient graph has vertex set $\mathbb F_p^2$ and steps
\[
\pm(1,0),\qquad\pm(0,1),\qquad\pm(1,1).
\]
A hamiltonian cycle in the quotient visits, for $j=0,\ldots,p-1$, the vertices
\[
(-j,j),\ (-j+1,j),\ \ldots,\ (-j+p-1,j),
\]
and then moves to $(-j-1,j+1)$.
Every step within a row is $(1,0)$, and the step to the next row is $(0,1)$ because $-j+p-1=-j-1$ in $\mathbb F_p$.
Each row contains all $p$ possible first coordinates, and the
last vertex $(0,p-1)$ is joined to $(0,0)$ by a $(0,1)$-step.
Thus this is a hamiltonian cycle.
List the corresponding orbits in this cyclic order as $B_0,\ldots,B_{p^2-1}$.
Explicitly, we have $B_{jp+r}=B_{-j+r,j}$ for $0\le j,r<p$.

Replace every edge of $X$ by its two opposite arcs, and keep
all arcs from $B_s$ to $B_{s+1}$, with indices modulo $p^2$.
Denote the resulting spanning digraph by $Y$.
We choose $(i,j,0,0)$ as the representative of $B_{i,j}$.
Under the identification $(d^ue^v,jp+r)\longmapsto(-j+r,j,u,v)$, equation~\eqref{eq:coset-p4-edges} gives
$Y\cong D(N;p^2;S_0,\ldots,S_{p^2-1})$,
where
\[
S_{jp+r}=
\begin{cases}
\langle d\rangle,&0\le r<p-1,\\
\langle e\rangle,&r=p-1,
\end{cases}
\qquad 0\le j<p.
\]
Since $N$ is abelian, we infer that
$S_0S_1\cdots S_{p^2-1}=\langle d\rangle\langle e\rangle=N$.
Therefore \Cref{clm:scc} shows that $Y$ is strongly connected, and \Cref{lem:cyclic-layers} gives a directed hamiltonian cycle in $Y$. 
Its underlying cycle is a hamiltonian cycle of $X$.

We next show that $X$ is not a Cayley graph.
Let $o=(0,0,0,0)$. For $t\in\mathbb F_p$, put
\[
\begin{aligned}
u_t&=(1,0,t,0),&
v_t&=(1,1,t,t),&
w_t&=(0,1,0,t),\\
u'_t&=(-1,0,t,0),&
v'_t&=(-1,-1,t,t+1),&
w'_t&=(0,-1,0,t).
\end{aligned}
\]
By \eqref{eq:coset-p4-edges}, these are precisely
the $6p$ neighbours of $o$, and their induced
adjacencies form the cycle
\[
u_t\edge v_t\edge w_t\edge u'_{-t}
\edge v'_{-t}\edge w'_{1-t}\edge u_{t-1},
\qquad t\in\mathbb F_p.
\]
These are all the possible adjacencies among the six families, as follows by comparing their first two coordinates. 
Since $t$ decreases by one after each six edges, $X[N_X(o)]\cong C_{6p}$.%

We put $A=\Aut(X)$ and claim that $A_o$ acts faithfully on $N_X(o)$. 
Suppose that an automorphism fixes $o$ and all its neighbours.
For every $u\in N_X(o)$, there is a common neighbour $v$ of $o$ and $u$. 
Thus $o,v$ are two adjacent fixed vertices of the cycle
induced by $N_X(u)$, so this entire cycle is fixed pointwise. 
Connectivity of $X$ then implies that the automorphism is the identity.

Thus we get  $A_o\le\Aut(C_{6p})$, where $|A_o|\mid 12p$.
Since $|V(X)|=p^4$ and $p\ge5$, the group $G$, of order $p^5$, is a Sylow $p$-subgroup of $A$.

Suppose that $X$ is a Cayley graph. 
Then $A$ contains a regular subgroup $R$ of order $p^4$. 
By Sylow's theorem, we may conjugate
$R$ so that $R\le G$. Thus $R$ has index
$p$ in $G$, and hence $\Phi(G)\le R$.
But $c=[a,b]\in\Phi(G)$, so $H=\langle c\rangle\le R$. 
Since $H$ fixes the vertex $H$, this contradicts the regularity of $R$. Therefore $X$
is not a Cayley graph.
\end{exa}

\section*{Data Availability}

No datasets were generated or analysed during the current study.

\section*{Declarations}

The authors declare no conflict of interest.

\end{document}